\documentclass[twoside,11pt]{article}
\usepackage{amsmath,amssymb,graphicx,amsthm}
\usepackage[latin1]{inputenc}
\usepackage{anysize}

\marginsize{2,5cm}{2,5cm}{2cm}{3cm}

\usepackage{tikz}

\theoremstyle{definition}
  \newtheorem{definition}{Definition}[section]

\theoremstyle{plain}
  \newtheorem{theorem}{Theorem}[section]

\newcommand{\re}{\mathbb{R}}

\title{Delta shock wave for a $3 \times 3$ hyperbolic system of conservation laws}
\author{Richard De la cruz\thanks{Universidad Pedagógica y Tecnológica de Colombia, Tunja, Colombia. {\tt richard.delacruz@uptc.edu.co}}\\
Juan Galvis\thanks{Universidad Nacional de Colombia, Bogotá, Colombia. {\tt jcgalvisa@unal.edu.co}}\\
Juan Carlos Juajibioy\thanks{Universidad Nacional de Colombia, Bogotá, Colombia. {\tt jcjuajibioyo@unal.edu.co}}\\
Leonardo Rendón\thanks{Universidad Nacional de Colombia, Bogotá, Colombia. {\tt lrendona@unal.edu.co}}}
\begin{document}
 \noindent
 \maketitle
 
 \begin{center}
  {\em XV International Conference on Hyperbolic Problems:\\ Theory, Numerics, Applications - HYP 2014}
 \end{center}
 
 \begin{abstract}
  We study the one-dimensional Riemann problem for a hyperbolic system of
three conservation laws of Temple class. The system  is a simplification of
a recently propose system of five conservations laws by Bouchut and Boyaval
that models viscoelastic fluids. An important issue is that the considered $3 \times 3$
system is such that every characteristic field is linearly degenerate. We study
the Riemann problem for this system and under suitable generalized Rankine-Hugoniot 
relation and entropy condition, both existence and uniqueness of particular delta-shock type solutions are established.
 \end{abstract}

 \section{Introduction}
 In this work, we show existence and uniqueness of delta shock wave for the Suliciu relaxation system \cite{BouchutBook, Suliciu1}
 \begin{equation}\label{sistema}
 \begin{cases}
  \rho_t+(\rho u)_x=0,\\
  (\rho u)_t+(\rho u^2+s^2 v)_x=0,\\
  (\rho v)_t+(\rho uv+u)_x=0,
 \end{cases}
 \end{equation}
 with $s=const. >0$. The concept of delta shock wave is a generalization  of an ordinry shock wave. 
This generalization was introdduced by Korchinski  in the year of 1977  in his PhD thesis \cite{Korchinski}. He considered the Riemann problem for system
 \begin{equation} \label{SisKor}
 \begin{cases}
  u_t+\left( \frac12 u^2 \right)_x=0,\\
  v_t+\left( \frac12 uv \right)_x=0.
 \end{cases}
 \end{equation}
Motivated by some numerical results, he constructed the unique Riemann solution using generalized delta functions to obtain singular shocks satisfying 
\eqref{SisKor} in the sense of distributions. After that, in 1994, Tan, Zhang and Zheng established in \cite{TanZhangZheng} the existence, uniqueness and stability of delta shock waves 
for a viscous perturbation of system \eqref{SisKor} where the first equation is replaced by 
$u_t+\left( u^2 \right)_x=\varepsilon u_{xx}$, with $\varepsilon>0$. Other works dealing with delta-shock waves are due to Ercole \cite{Ercole} who in 2000 obtained a delta shock solution as a limit of smooth solutions by the vanishing viscosity method.  Sheng and Zhang \cite{ShengZhang} discussed the Riemann problem for pressureless gas system.
 In 2005, Brenier \cite{YBrenier} considered the Riemann problem for the Chaplygin gas system. Other works for the Chaplygin gas system can be found in \cite{ChengYang, Guo, KongWei}.
 Some more recent works on delta shocks for general hyperbolic conservation laws are due to Danilov and Mitrovic \cite{DanilovMitrovic, Danilov-Shelkovich}, where they described delta shock wave generation from continuous initial data by using smooth approximations in the weak sense. 
 { The works mentioned above on the delta shock solutions are for particular
cases of $2\times 2$ systems. For particular $3\times 3$ non-strictly hyperbolic systems with repeated eigenvalues, we mention the works of Panov and Shelkovich \cite{Panov-Shelkovich}, Shelkovich \cite{Shelkovich1} and Cheng \cite{HCheng}. For a strictly hyperbolic system, see \cite{DelacruzG} for  partial results of delta shock solution.}
For new developments of delta shock waves, its
applications and historical notes is \cite{avancesDelta}.\\

Coming back to system \eqref{sistema}, we can deduce that the associated eigenvalues are given by,
$\lambda_1 = u - \frac{s}{\rho}$,  $\lambda_2 = u$ and $\lambda_3 = u + \frac{s}{\rho}$ 
where the corresponding Riemann invariants are
$R_1 = s^2v-su$, $R_2 = v + \frac{1}{\rho}$ and $R_3 = s^2v+su$. From the expressions for the eigenvalues and the Riemann invariants we obtain
$\lambda_1 = \frac{R_3}{s}-sR_2$, $\lambda_2 = \frac{1}{2s}(R_3-R_1)$ and $\lambda_3 = sR_2-\frac{R_1}{s}$.
We can see that system \eqref{sistema} is linearly degenerate.\\
 
 Recently, Lu et al. \cite{Lu} showed existence of solutions for the Cauchy problem associated to the Suliciu relaxation system \eqref{sistema} with bounded initial data 
 \begin{equation} \label{datoincial}
  \begin{aligned}
   & (\rho(0,x),u(0,x),v(0,x))=(\rho_0(x),u_0(x),v_0(x)),  \quad \rho_0(x) \ge 0, \quad x \in \re,
  \end{aligned}
 \end{equation}
subject to the following conditions:
\begin{enumerate}
 \item[(H1)] The functions $\rho_0$, $u_0$ and $v_0$ satisfy
 \begin{align*}
  c_1 \le u_0(x)-sv_0(x) \le c_2, \qquad c_3 \le u_0(x)+sv_0(x) \le c_4, \qquad
  v_0(x)+\frac{1}{\varepsilon+ \rho_0(x)}>c_5, 
 \end{align*}
where $\varepsilon>0$ is a small and the constants $c_i$ , $i = 1, \dots , 5$, satisfy
$ c_5-\frac{c_4-c_1}{2s}>0$.

\item[(H2)] The total variations of $u_0(x) - sv_0(x)$ and $u_0(x) +s v_0(x)$ are bounded.
\end{enumerate}
The existence result for the Cauchy problem \eqref{sistema}--\eqref{datoincial} includes solutions in vacuum regions \cite[Theorem 1]{Lu}. When $\rho_0(x) \ge \underline{\rho} >0$, the condition H1 becomes:
\begin{enumerate}
 \item[(H1)] The functions $\rho_0$, $u_0$ and $v_0$ satisfy
$ c_1 \le u_0(x)-sv_0(x) \le c_2$, $c_3 \le u_0(x)+sv_0(x) \le c_4$, 
and $ v_0(x)+\frac{1}{\rho_0(x)}>c_5$, 
where $c_i$ , $i = 1, \dots , 5$, are  constants satisfying
$ c_5-\frac{c_4-c_1}{2s}>0.$
\end{enumerate}
In \cite{DelacruzG}, we showed uniqueness of solutions for the Cauchy problem associated to the Suliciu relaxation system \eqref{sistema}--\eqref{datoincial} with $\rho_0(x) \ge \underline{\rho} >0$. Moreover, 
Theorem 2 in \cite{DelacruzG} ensures that when we consider $v_0(x)=-\frac{1}{\rho_0(x)}$, problem 
\eqref{sistema} is reduced to following Chaplygin gas system
\begin{equation*}
 \begin{cases}
  \rho_t+(\rho u)_x=0,\\
  (\rho u)_t+(\rho u^2-\frac{s^2}{\rho})_x=0.
 \end{cases}
\end{equation*}
The classical Riemann problem associated to the Suliciu relaxation system \eqref{sistema}--\eqref{datoIniRiemann} has been extensively studied, for instance in \cite{BouchutBook, Carbou, Chalons-Coquel, DelacruzG}. In \cite{DelacruzG1}, it is shown uniqueness for the generalized Riemann problem for the Suliciu relaxation system.

In the present paper, we show the existence and uniqueness of solutions for the Riemann problem associated with the Suliciu relaxation system
with initial data 
\begin{equation}\label{datoIniRiemann}
 (\rho_0(x),u_0(x),v_0(x))=\begin{cases}
                                                               (\rho_l,u_l,v_l), & \text{if } x<0,\\
                                                               (\rho_r,u_r,v_r), & \text{if } x>0,
                                                           \end{cases}
\end{equation}
in which the left and right constant states $(\rho_l,u_l,v_l)$ and $(\rho_r,u_r,v_r)$, with $\rho_l, \rho_r >0$, satisfy the conditions H2 
 and $\lambda_1(\rho_l,u_l,v_l) \ge \lambda_3(\rho_r,u_r,v_r)$, i.e., they do not satisfy condition H1 globally.
Under this assumptions  we have the following situations:
\begin{enumerate}
 \item $(\rho_l,u_l,v_l)$ and $(\rho_r,u_r,v_r)$ satisfy locally the condition H1, i.e.,
 \begin{itemize}
  \item $(\rho_l,u_l,v_l)$ satisfy
  \begin{equation*}
   \begin{aligned}
    \alpha_1 \le u_l-sv_l \le \alpha_2, \qquad \alpha_3 \le u_l-sv_l \le \alpha_4 \quad
        \text{ and } \quad v_l+\frac{1}{\rho_l} > \alpha_5
   \end{aligned}
  \end{equation*}
where $\alpha_i$, $i=1,\dots,5$, are suitable constants satisfying $\alpha_5-\frac{\alpha_4-\alpha_1}{2s}>0$,
  \item $(\rho_r,u_r,v_r)$ satisfy
  \begin{equation*}
   \begin{aligned}
    \beta_1 \le u_r-sv_r \le \beta_2, \qquad \beta_3 \le u_r-sv_r \le \beta_4 \quad
        \text{ and } \quad v_r+\frac{1}{\rho_r} > \beta_5
   \end{aligned}
  \end{equation*}
where $\beta_i$, $i=1,\dots,5$, are suitable constants satisfying $\beta_5-\frac{\beta_4-\beta_1}{2s}>0$.
\item Let $c_1=\min \{ \alpha_1,\beta_1 \}$ and $c_4=\max \{ \alpha_1,\beta_1 \}$. Then there are $c_i$, $i=1,\dots,4$ such that
\begin{equation*}
  c_1 \le \begin{Bmatrix}
           u_l-sv_l \\ u_r-sv_r
          \end{Bmatrix}
\le c_2
\quad \text{ and } \quad
  c_3 \le \begin{Bmatrix}
           u_l+sv_l \\ u_r+sv_r
          \end{Bmatrix}
\le c_4,
\end{equation*}
but is not possible to find a constant $c_5$ such that $c_5-\frac{c_4-c_1}{2s}>0$.
 \end{itemize}
 \item Only $(\rho_r,u_r,v_r)$ satisfy locally the condition H1, i.e., $(\rho_l,u_l,v_l)$ satisfy
  \begin{equation*}
   \begin{aligned}
    \beta_1 \le u_l-sv_l \le \beta_2, \qquad \beta_3 \le u_l-sv_l \le \beta_4 \quad
        \text{ and } \quad v_l+\frac{1}{\rho_l} > \beta_5
   \end{aligned}
  \end{equation*}
where $\beta_i$, $i=1,\dots,5$, are suitable constants satisfying $\beta_5-\frac{\beta_4-\beta_1}{2s} \le 0$.
 \item Only $(\rho_l,u_l,v_l)$ satisfy locally the condition H1, i.e., $(\rho_r,u_r,v_r)$ satisfy
  \begin{equation*}
   \begin{aligned}
    \beta_1 \le u_r-sv_r \le \beta_2, \qquad \beta_3 \le u_r-sv_r \le \beta_4 \quad
        \text{ and } \quad v_r+\frac{1}{\rho_r} > \beta_5
   \end{aligned}
  \end{equation*}
where $\beta_i$, $i=1,\dots,5$, are suitable constants satisfying $\beta_5-\frac{\beta_4-\beta_1}{2s} \le 0$.
 \item Neither $(\rho_l,u_l,v_l)$ nor $(\rho_r,u_r,v_r)$ satisfy the local condition H1.
\end{enumerate}
De la cruz et al. \cite{DelacruzG} studied only the first case. The idea used by the authors in \cite{DelacruzG} can be extended to the other cases. 

We considered all the cases for which existence and uniqueness of solutions 
can be guarantied for the Riemann problem associated to the Suliciu system.
In \cite{DelacruzG} it is shown that if the initial data \eqref{datoincial}, with $\rho_0(x) \ge \underline{\rho}>0$, satisfies H1  and  H2, then, the Riemann problem has a unique solution. It is also shown that  H1 and H2 imply the 
Lax shock condition $\lambda_1(\rho_l,u_l,v_l)<\lambda_3(\rho_l,u_l,v_l)$. 
Therefore, the Riemann problem for the Suliciu relaxation system has a classical 
solution in the region $\Gamma_1=\{ (\rho,u,v) \, : \, \rho>0 \text{ and } \lambda_1(\rho_l,u_l,v_l)<\lambda_3(\rho,u,v) \}$. 
In this paper we show that the Suliciu relaxation system presents nonclassical solutions
in the sense of delta shocks  in the region $$\Gamma_2=\{ (\rho,u,v) \, : \, \rho>0, \lambda_1(\rho_l,u_l,v_l) \ge \lambda_3(\rho,u,v) \text{ and } (u_l-u)^2 \ge s^2(v_l-v)(1/\rho-1/\rho_l)\}.$$  

It remains to study if we can obtain some kind of solution in 
the region $(0,\infty)\times \re^2-(\Gamma_1\cup \Gamma_2)$. This object of our 
current research.

\section{Delta shock solutions}
Denote by $BM(\re)$ the space of bounded Borel measures on $\re$. the definition of a measure solution of Suliciu relaxation system in $BM(\re)$ can be given as follows.
\renewcommand{\labelenumi}{\alph{enumi}$)$ }
\begin{definition} \label{measureSol}
 A triple $(\rho,u,v)$ constitutes a \emph{measure solution} to the Suliciu relaxation system, if it holds that
 \begin{enumerate}
  \item $\rho \in L^\infty((0,\infty),BM(\re)) \cap C((0,\infty),H^{-s}(\re))$,
  \item $u \in L^\infty((0,\infty),L^\infty(\re)) \cap C((0,\infty),H^{-s}(\re))$,
  \item $v \in L_{loc}^\infty((0,\infty),L_{loc}^\infty(\re)) \cap C((0,\infty),H^{-s}(\re))$, \ \ $s>0$,
  \item $u$ and $v$ are measurable with respect to $\rho$ at almost for all $t\in(0,\infty)$,
 \end{enumerate}
and
\begin{equation} \label{MS1A}
\begin{cases}
 I_1=\int_0^\infty \int_{\re} (\phi_t+u\phi_x) \, d\rho dt =0,\\
 I_2=\int_0^\infty \int_{\re} u(\phi_t+u\phi_x) \, d\rho dt + \int_0^\infty \int_{\re} s^2 v \phi_x \, dx dt =0, \\
 I_3=\int_0^\infty \int_{\re} v(\phi_t+u\phi_x) \, d\rho dt + \int_0^\infty \int_{\re} u \phi_x \, dx dt =0,
\end{cases}
\end{equation}
for all test function $\phi \in C_0^\infty(\re^+ \times \re)$.
\end{definition}

Recall that a two-dimensional weighted delta function $w(s)\delta_L$ supported on a smooth curve $L$ parameterized as $t=t(s)$, $x=x(s)$ $(c \le s \le d)$ is defined by
$  \langle w(s)\delta_L, \phi(t,x) \rangle = \int_c^d w(s)\phi(t(s),x(s)) \, ds
$ for all $\phi \in C_0^\infty(\re^2)$. We next present the definition of 
a delta shock wave.

\begin{definition}
 A triple distribution $(\rho,u,v)$ is called a \emph{delta shock wave} if it is represented in the form
 \begin{equation} \label{deltashock}
  (\rho,u,v)(t,x)=\begin{cases}
                   (\rho_l,u_l,v_l)(t,x), & x<x(t),\\
                   (w(t)\delta(x-x(t)),u_\delta(t),g(t)), & x=x(t),\\
                   (\rho_r,u_r,v_r)(t,x), & x>x(t),
                  \end{cases}
 \end{equation}
and satisfies Definition \ref{measureSol}, where $(\rho_l,u_l,v_l)(t,x)$ and $(\rho_r,u_r,v_r)(t,x)$ are piecewise smooth bounded solutions of 
the Suliciu relaxation system \eqref{sistema}.
\end{definition}

We set $\frac{dx}{dt}=u_\delta(t)$ since the concentration in $\rho$ need to travel at the speed of discontinuity. Hence, we say that a delta shock wave \eqref{deltashock} is a measure solution to the Suliciu relaxation system \eqref{sistema} if and only if the
following relation holds,
\begin{equation} \label{RHgener}
 \begin{cases}
  \frac{dx(t)}{dt}=u_\delta(t),\\
  \frac{dw(t)}{dt}=-[\rho]u_\delta(t)+[\rho u],\\
  \frac{dw(t)u_\delta(t)}{dt}=-[\rho u]u_\delta(t)+[\rho u^2+s^2v],\\
  \frac{dw(t)g(t)}{dt}=-[\rho v]u_\delta(t)+[\rho uv+u].
 \end{cases}
\end{equation}
In fact, for any test function $\phi \in C_0^\infty (\re^+ \times \re)$,
from \eqref{MS1A}, we obtain
\begin{equation*}
 \begin{aligned}
  I_1 &=\int_0^\infty \int_{\re} (\phi_t+u\phi_x) \, d\rho dt =\int_0^\infty \left\{ -u_\delta(t)[\rho]+[\rho u]-\frac{dw(t)}{dt} \right\} \phi \, dt,\\
 I_2 &
 = \int_0^\infty \left\{ -u_\delta(t)[\rho u]+[\rho u^2+s^2v]-\frac{dw(t)u_\delta(t)}{dt} \right\} \, dt, \quad \text{ and } \\
 I_3 
 &=\int_0^\infty \left\{ -u_\delta(t)[\rho v]+[\rho uv+u]-\frac{dw(t)g(t)}{dt} \right\} \phi \, dt.
 \end{aligned}
\end{equation*}
Relations \eqref{RHgener} are called the {\em generalized Rankine-Hugoniot condition}. It  reflects the exact relationship among the
limit states on two sides of the discontinuity, the weight,
propagation speed and the location of the discontinuity.
In addition, to guarantee uniqueness, the delta shock wave should satisfy 
the admissibility (entropy) condition
$ \lambda_3(\rho_r,u_r,v_r) \le u_\delta(t) \le \lambda_1(\rho_l,u_l,v_l).$\\

Now, the generalized Rankine-Hugoniot condition \eqref{RHgener} is applied to the Riemann problem \eqref{sistema}--\eqref{datoIniRiemann} with left and right
 constant states $U_-=(\rho_-,u_-,v_-)$ and
$U_+=(\rho_+,u_+,v_+)$, respectively, satisfying the condition H2, the fact $\lambda_3(\rho_+,u_+,v_+) \le \lambda_1(\rho_-,u_-,v_-)$ and
\begin{equation} \label{CondicionDeltaSol}
\begin{aligned}
 (u_--u_+)^2 \ge \frac{s^2}{\rho_+}(v_--v_+) - \frac{s^2}{\rho_-}(v_--v_+).
\end{aligned}
\end{equation}

Therefore, the Riemann problem \eqref{sistema}--\eqref{datoIniRiemann} is reduced to solving \eqref{RHgener} with initial data
\begin{equation} \label{datoRHgener}
 t=0, \quad x(0)=0, w(0)=0, g(0)=0,
\end{equation}
under entropy condition
\begin{equation} \label{condEntroGen}
 u_++\frac{s}{\rho_+} \le u_\delta(t) \le u_--\frac{s}{\rho_-}.
\end{equation}
From \eqref{RHgener} and \eqref{datoRHgener}, it follows that
\begin{equation}\label{DFREW}
 \begin{aligned}
  w(t) &= -[\rho]x(t)+[\rho u]t, \\
  w(t)u_\delta(t) &= -[\rho u]x(t)+[\rho u^2+s^2v]t, \text{  and } \\
  w(t)g(t) &= -[\rho v]x(t)+[\rho uv+u]t.
 \end{aligned}
\end{equation}
Multiplying the first equation in \eqref{DFREW} by $u_\delta(t)$ and then subtracting it from the second one, we obtain that
\begin{equation*}
 [\rho]x(t)u_\delta(t)-[\rho u]u_\delta(t)t-[\rho u]x(t)+[\rho u^2+s^2v]t=0,
\end{equation*}
that is,
\begin{equation*}
 \frac{d}{dt} \left( \frac{[\rho]}{2}x^2(t)-[\rho u]x(t)t+\frac{[\rho u^2+s^2v]}{2}t^2 \right) =0,
\end{equation*}
which is equivalent to
\begin{equation} \label{cuadratica1}
 [\rho]x^2(t)-2[\rho u]x(t)t+[\rho u^2+s^2v]t^2 =0.
\end{equation}
From \eqref{cuadratica1}, one can find $u_\delta(t) := u_\delta$ is a constant and $x(t) = u_\delta t$. Then, \eqref{cuadratica1} can
be rewritten
\begin{equation} \label{cuadratica2}
 [\rho]u_\delta^2-2[\rho u]u_\delta+[\rho u^2+s^2v]=0.
\end{equation}

When $[\rho]=\rho_--\rho_+=0$, the situation is very simple and one can easily calculate the solution
\begin{equation} \label{SolDeltaLin}
 \begin{cases}
  u_\delta=\frac{u_-+u_+}{2}+s^2\frac{[v]}{2\rho_-[u]}, \\
  x(t)= u_\delta t, \\
  w(t)=\rho_-(u_--u_+)t, \\
  g(t)= \frac{[\rho u v + u]-u_\delta}{[\rho u]},
 \end{cases}
\end{equation}
which obviously satisfies the entropy condition \eqref{condEntroGen}.
From condition \eqref{CondicionDeltaSol},
\begin{align*}
 s^2 \frac{[v]}{\rho_-} \le \frac12 (\lambda_1(U_-)-\lambda_3(U_+))^2 < \frac12 [u] (\lambda_1(U_-)-\lambda_3(U_+))
\end{align*}
and
\begin{align*}
 u_\delta-\left(u_--\frac{s}{\rho_-} \right) 
 =\frac12 \left(  \left( u_++\frac{s}{\rho_-} \right) - \left( u_--\frac{s}{\rho_-} \right) + s^2 \frac{[v]}{\rho_-[u]} \right) \le 0.
\end{align*}
Similarly, we can deduce that
\begin{align*}
 u_\delta-\left(u_++\frac{s}{\rho_-} \right) 
 =\frac12 \left(  \left( u_--\frac{s}{\rho_-} \right) - \left( u_++\frac{s}{\rho_-} \right) + s^2 \frac{[v]}{\rho_-[u]} \right) \ge 0,
\end{align*}
because
\begin{align*}
 -s^2 \frac{[v]}{\rho_-} \le \frac12 (\lambda_1(U_-)-\lambda_3(U_+))^2 < \frac12 [u] (\lambda_1(U_-)-\lambda_3(U_+)).
\end{align*}
When $[\rho]=\rho_--\rho_+ \neq 0$, the discriminant of the quadratic equation \eqref{cuadratica2} is
\begin{equation*}
 \Delta = 4[\rho u]^2-4[\rho][\rho u^2+s^2v]=\rho_-\rho_+[u]^2-s^2[\rho][v]>0
\end{equation*}
and then we can find
\begin{equation} \label{SolDelta1}
(S^\pm)=
 \begin{cases}
  u_\delta=\frac{[\rho u]\mp \sqrt{\Delta/4}}{[\rho]},\\
  x(t)=\frac{[\rho u]\mp \sqrt{\Delta/4}}{[\rho]} t,\\
  w(t)= \pm \sqrt{\Delta/4} t,\\
  g(t)= \frac{-[\rho u][\rho v]\pm [\rho v]\sqrt{\Delta/4}+[\rho][\rho u v +u]}{[\rho]\sqrt{\Delta/4}}t.
 \end{cases}
\end{equation}

Next, with the help of the entropy condition \eqref{condEntroGen}, we will choose the admissible solution from \eqref{SolDelta1}. 
Observe that by the entropy condition and since the system is strictly hyperbolic, we have that
$$ u_+-\frac{s}{\rho_+}<u_+<u_++\frac{s}{\rho_+} \le u_--\frac{s}{\rho_-} <u_-<u_-+\frac{s}{\rho_-}. $$
Observe that,
\begin{align*}
 -[\rho]\lambda_1(\rho_-,u_-,v_-)+[\rho u]=\rho_+ \left( \left(u_--\frac{s}{\rho_-} \right) - \left(u_+-\frac{s}{\rho_+} \right) \right) >0,\\
 -[\rho]\lambda_3(\rho_+,u_+,v_+)+[\rho u]=\rho_- \left( \left(u_-+\frac{s}{\rho_-} \right) - \left(u_++\frac{s}{\rho_+} \right) \right) >0,
\end{align*}
\begin{align*}
 [\rho](\lambda_1(\rho_-,u_-,v_-))^2-2[\rho u]\lambda_1(\rho_-,u_-,v_-)+[\rho u^2]+s^2[v] = \\
 -\rho_+\left( u_--u_+-\frac{s}{\rho_-} \right)^2+\frac{s^2}{\rho_-}+s^2[v] \le 0,
\end{align*}
\begin{align*}
 [\rho](\lambda_3(\rho_+,u_+,v_+))^2-2[\rho u]\lambda_3(\rho_+,u_+,v_+)+[\rho u^2]+s^2[v] = \\
 \rho_-\left( u_--u_+-\frac{s}{\rho_+} \right)^2-\frac{s^2}{\rho_+}+s^2[v] \ge 0,
\end{align*}
then, for the solution given in \eqref{SolDelta1} by $(S^+)$, we have
\begin{align*}
 u_\delta-  \lambda_1(\rho_-,u_-,v_-) 
 \le 0
\qquad \text{ and } \qquad
 u_\delta-  \lambda_3(\rho_+,u_+,v_+) 
 \ge 0,
\end{align*}
which imply that the entropy condition \eqref{condEntroGen} is valid.
When $\lambda_1(\rho_-,u_-,v_-)=\lambda_3(\rho_+,u_+,v_+)$, we have trivially that
$\lambda_1(\rho_-,u_-,v_-)=u_\delta=\lambda_3(\rho_+,u_+,v_+)$.

Now, for the solution $(S^-)$, 
when $\rho_-<\rho_+$ we have
\begin{equation*}
 \begin{aligned}
  u_\delta &-\lambda_3(\rho_+,u_+,v_+)
  =\frac{\rho_-(\lambda_3(U_-)-\lambda_3(U_+))+\sqrt{[\rho u]^2-[\rho][\rho u^2+s^2v]}}{[\rho]} <0,
 \end{aligned}
\end{equation*}
and when $\rho_->\rho_+$, that
\begin{equation*}
 \begin{aligned}
  u_\delta &-\lambda_1(\rho_-,u_-,v_-)
  =\frac{\rho_+(\lambda_1(U_-)-\lambda_1(U_+))+\sqrt{[\rho u]^2-[\rho][\rho u^2+s^2v]}}{[\rho]} >0.
 \end{aligned}
\end{equation*}
showing that the solution $(S^-)$ in \eqref{SolDelta1} 
does not satisfy the entropy
condition \eqref{condEntroGen}.\\
Thus we have proved the following result.
\begin{theorem}
 Given left and right constant states $(\rho_l,u_l,v_l)$ and $(\rho_r,u_r,v_r)$,
respectively, such that satisfy the condition H2, $\lambda_1(\rho_l,u_l,v_l) \ge \lambda_3(\rho_r,u_r,v_r)$
and \eqref{CondicionDeltaSol}, that is,
\begin{equation*}
 (u_l-u_r)^2 \ge s^2(v_l-v_r)/\rho_r - s^2(v_l-v_r)/\rho_l.
\end{equation*}
Then, the Riemann problem \eqref{sistema}--\eqref{datoIniRiemann} admits a unique entropy solution in the sense of measures.
This solution is of the form
\begin{equation*}
 (\rho,u,v)(t,x)=\begin{cases}
                  (\rho_l,u_l,v_l), & \text{if } x<u_\delta t, \\
                  (w(t)\delta(x-u_\delta t),u_\delta,g(t)), & \text{if } x=u_\delta t, \\
                  (\rho_r,u_r,v_r), & \text{if } x>u_\delta t,
                 \end{cases}
\end{equation*}
where $u_\delta$, $w(t)$ and $g(t)$ are show in \eqref{SolDeltaLin} for $[\rho]=0$ or $(S^+)$ in \eqref{SolDelta1} for $[\rho] \neq 0$.
\end{theorem}

%

\section{Numerical illustration}
\begin{figure}[h]
  \centering
    \includegraphics[scale=0.3]{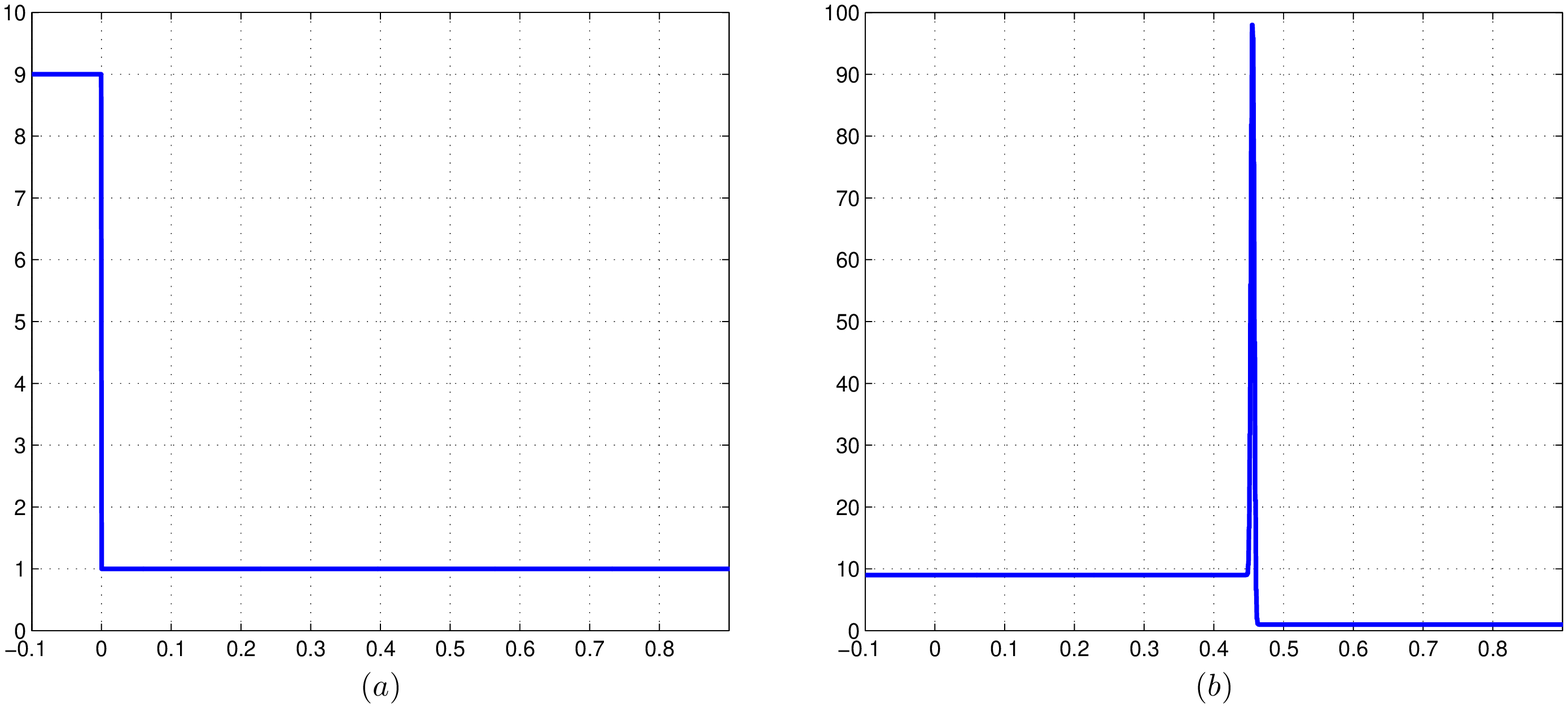}    
    \includegraphics[scale=0.3]{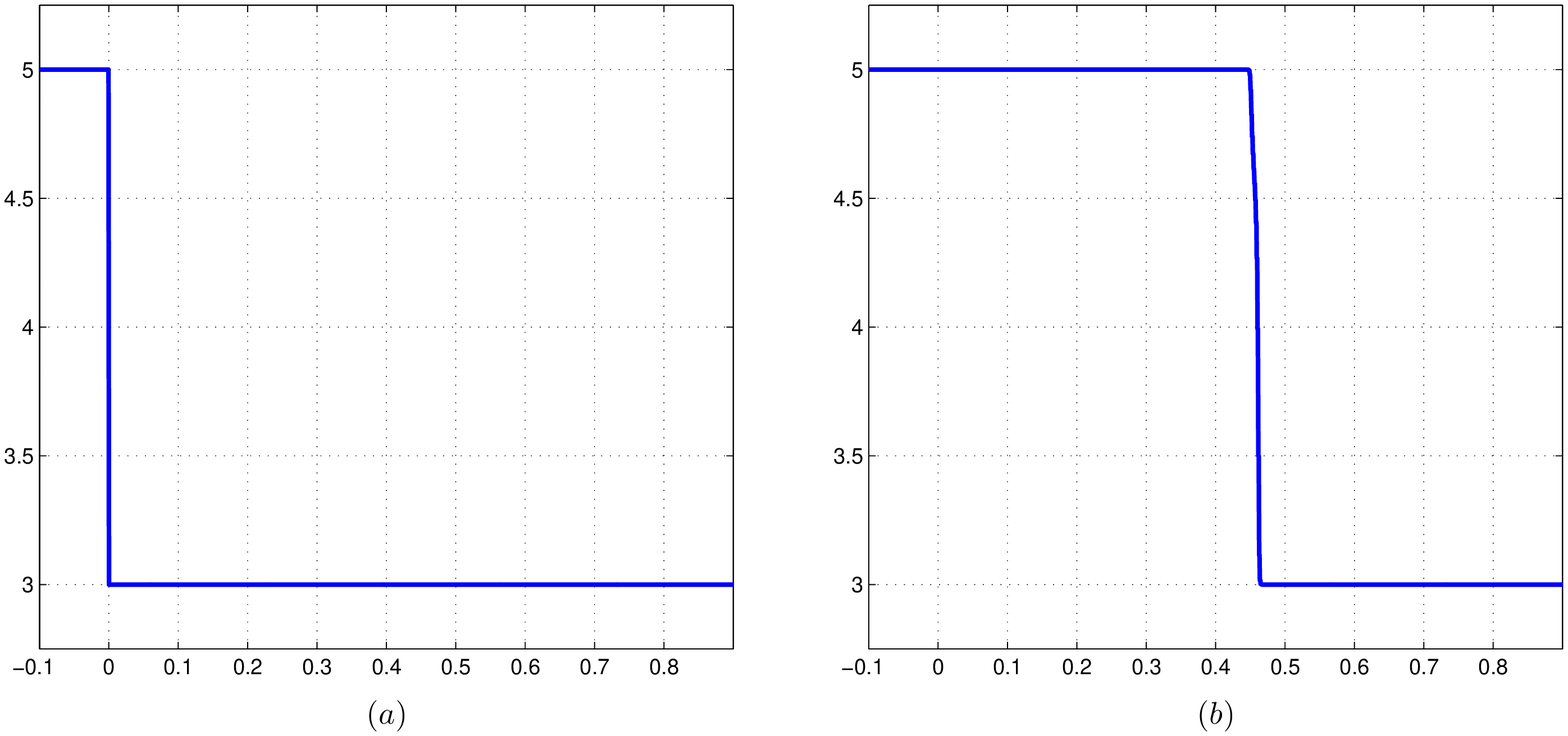}
    \includegraphics[scale=0.3]{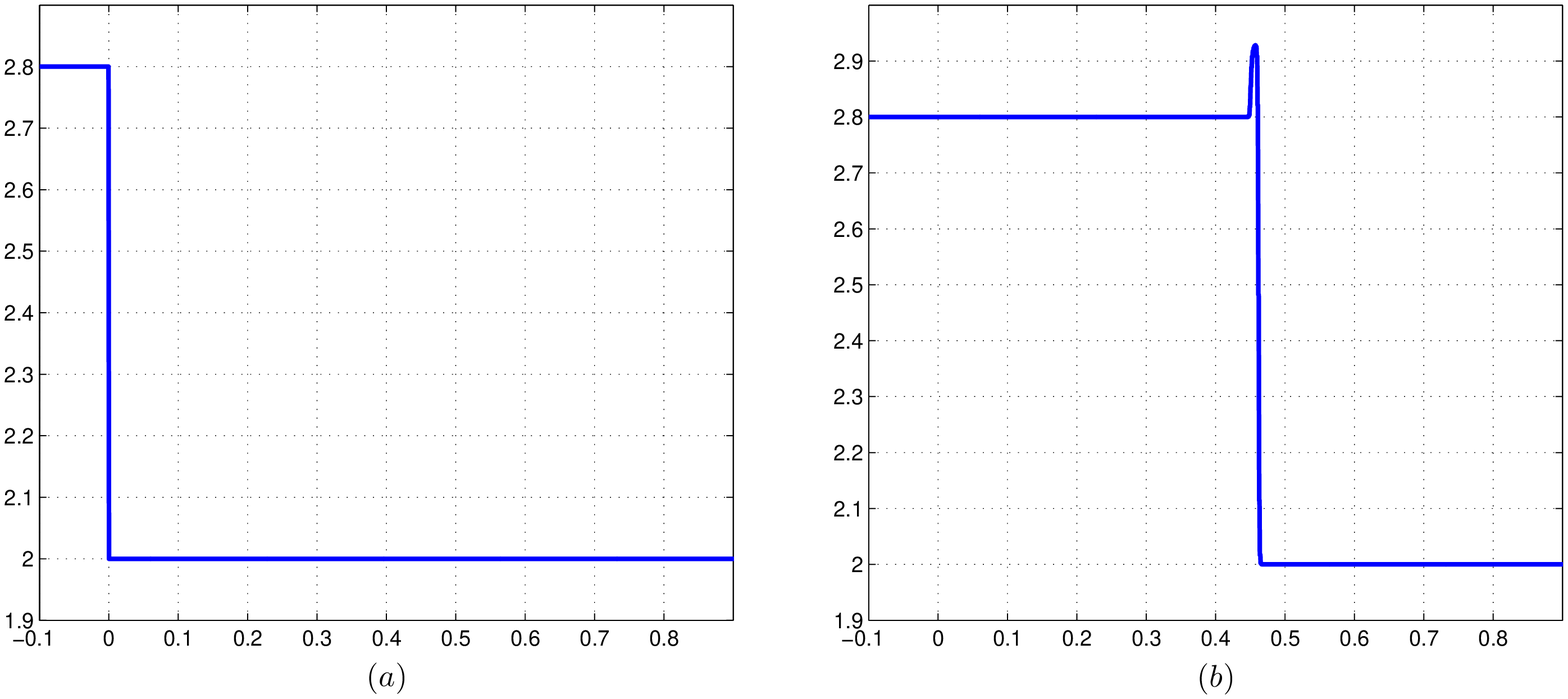}
  \caption{\emph{Numerical solution}. Top: (a) The initial data $\rho_0(x) = \rho(0, x)$. \,
(b) Numerical solution of $\rho$ at time $t = 0.1$. 
Middle: (a) The initial data $u_0(x) = u(0, x)$. \,
(b) Numerical solution of $u$ at time $t = 0.1$. Bottom:
(a) The initial data $v_0(x) = v(0, x)$. \,
(b) Numerical solution of $v$ at time $t = 0.1$
}
  \label{delta1}
\end{figure}

In this section, we show numerical evidence of delta shock solution for the Suliciu relaxation system using the Lax-Friedrichs method.
In the numerical test, with $s=1$, we consider the initial data given by
\begin{equation*}
 (\rho_0,u_0,v_0)(x)=\begin{cases}
                      (9,5,14/5), &\text{ if } x<0,\\
                      (1,3,2), &\text{ if } x>0,
                     \end{cases}
\end{equation*}
the spatial discretization parameter for $N = 1780$ points and a constant
$CFL = 0.1969889$.\\ The numerical results at final time $t = 0.1$ is  displayed in
 Figure \ref{delta1}.
The exact solution at time $t$ is
\begin{equation*}
 (\rho,u,v)(t,x)=\begin{cases}
                      (9,5,14/5), &\text{ if } x<u_\delta t,\\
                      (At\delta(x-u_\delta t),u_\delta, Bt), &\text{ if } x=u_\delta t,\\
                      (1,3,2), &\text{ if } x>u_\delta t,
                     \end{cases}
\end{equation*}
with $u_\delta=\frac{21\sqrt{5}-\sqrt{37}}{4\sqrt{5}}$, $A=\frac{2\sqrt{37}}{\sqrt{5}}$ and $B=\frac{\sqrt{5}+29\sqrt{37}}{10\sqrt{37}}$.\\

%
%

\section{Conclusions}
In this paper we obtain delta-shock waves for the Riemann problem associated 
to the Soliciu relaxation system. This work complements and extend the 
partial results concerning delta-shock waves for the Riemann problem for the Suliciu relaxation system that were obtained in \cite{DelacruzG}.
{We mention that in \cite{BouchutBook, Carbou} the classical Riemann problem for the Suliciu relaxation system
was solved.  In our analysis we consider all possible cases for the existence and uniqueness of delta-shock solution for the Suliciu relaxation system and by the Lax-Friedrichs method we show numerical evidence of such solutions. 
We identified particular regions where we can present classical solutions and delta 
shock waves.
However, we observed that there are regions in which the Riemann problem for this system does not have classical solution neither delta shock waves, which encourages 
us to seek solutions in some different sense to be able to analyze this missing 
regions.}

\bibliography{bibliogra}
\bibliographystyle{abbrv}

 \end{document}